\documentclass[11pt,reqno]{article}

\usepackage{amsmath,amsthm,amsfonts,amssymb}
\usepackage{t1enc}
\newtheorem{thm}{Theorem}[section]

\theoremstyle{definition}

\theoremstyle{remark}
\newtheorem{rem}{Remark}[section]

\def\orig{e}

\newcommand{\eps}{\varepsilon}

\newcommand{\T}{\mathcal{T}}

\newcommand{\R}{\mathbb{R}}
\newcommand{\Z}{\mathbb{Z}}
\newcommand{\N}{\mathbb{N}}
\newcommand{\1}{\textbf{1}}
\newcommand{\ro}{\textbf{r}}
\renewcommand{\P}{\mathbb{P}}
\newcommand{\PP}{\textbf{P}}
\newcommand{\EE}{\textbf{E}}
\newcommand{\E}{\mathbb{E}}

\newcommand{\ZZ}{\mathcal{Z}}

\newcommand{\eval}[2][\right]{\relax
  \ifx#1\right\relax \left.\fi#2#1\rvert}

\title{Dynamical sensitivity of recurrence and transience of branching random walks}
\author{Sebastian M\"{u}ller}
\begin{document}
\maketitle  {\abstract Consider a sequence of i.i.d. random variables $X_n$ where each random variable is
refreshed independently according to a Poisson clock. At any fixed time $t$ the law of the sequence is the same
as for the sequence at time $0$ but at random times  almost sure properties of the sequence may be violated. If
there are such \emph{exceptional times} we say that the property is \emph{dynamically sensitive}, otherwise we
call it \emph{dynamically stable}. In this note we consider branching random walks on Cayley graphs and prove
that recurrence and transience are dynamically stable in the sub-and supercritical regime. While the critical
case is left open in general we prove dynamical stability for a specific class of Cayley graphs. Our proof
combines techniques from the theory of branching random walks with those of dynamical percolation.}

\section{Introduction}
In Benjamini et al. \cite{MR1959784} several properties of i.i.d. sequences are studied in the dynamical point
of view.  In particular, it is proven that transience of the simple random walk on the lattice $\Z^d$ is
dynamically stable for $d\ge 5,$ and dynamically sensitive for $d=3,4$. While recurrence is dynamically stable
for $d=1$, see \cite{MR1959784}, it is dynamically sensitive in dimension $d=2$, see Hoffman \cite{MR2235173}.
Khoshnevisan studied in \cite{Khoshnevisan05} and \cite{Khoshnevisan06} other properties of dynamical random
walks. We also refer to a recent survey \cite{steif09} on dynamical percolation. In this note we define
dynamical branching random walks on Cayley graphs and study dynamical sensitivity of  recurrence and transience,
see Theorem \ref{thm:no_except}.

\subsection{Random Walks}
Let us first collect the necessary notations for random walks on groups; for more details we refer to
\cite{woess}. Let $G$ be a finitely generated group with group identity $\orig,$ the group operations are
written multiplicatively (unless $G$ is abelian). Let $q$ be a probability measure on a finite generating set of
$G$. The random walk on $G$ with law $q$ is the Markov chain with state space $G$ and transition probabilities
$p(x,y)=q(x^{-1}y)$ for $x,y\in G$. Equivalently, the process can be described on the product space $(G,q)^\N$:
the $n$-th projections $X_n$ of $G^\N$ onto $G$ constitute a sequence of independent $G$-valued random variables
with common distribution $q$. Hence the random walk starting in $x\in G$ can be described as
\[S_n=x X_1\cdots X_n,\quad n\geq 0.\]
If not mentioned otherwise the random walk starts in the group identity $\orig$. Let $P$ denote the transition
kernel of the random walk and $p^{(n)}(x,y)=\P(S_n=y|S_0=x)$ be the probability to go from $x$ to $y$ in $n$
steps. We will assume the random walk to be irreducible, i.e., for all $x,y$ there exists some $k$ such that
$p^{(k)}(x,y)>0.$ Denote $G(x,y|z)=\sum_{n=0}^\infty p^{(n)}(x,y) z^n$ the corresponding generating functions.
The inverse of the convergence radius  of $G(x,y|z)$ is denoted by $\rho(P).$ The spectral radius $\rho(P)$ is
also given as
\begin{equation}\label{eq:rho:prop}
\rho(P)=\limsup_{n\to\infty} \left(p^{(n)}(x,y) \right),
\end{equation}
where the $\limsup$ is, due to the irreducibility of the random walk,  independent of $x$ and $y$.

\subsection{Branching Random Walks}
We use the interpretation of tree-indexed random walks, compare with \cite{benjamini:94}, to define the
branching random walk (BRW). Let $\T$ be a tree with root $\ro$. For a vertex $v$ of $\T$ let $|v|$ be the
(graph) distance from $v$ to the root $\ro$. Let $\T_n=\{v: |v|=n\}$ be the $n$-th level of the tree $\T$. We
label the edges of $\T$ with i.i.d. random variables with distribution $q$. The random variable $X_v$ is the
label of the edge $\langle v^-,v\rangle$ where $v^-$ is the unique predecessor of $v$, i.e., $|v^-|=|v|-1$.
Define $S_v=\orig \cdot \prod_{i=1}^{|v|} X_{v_i}$ where $\langle v_0=\ro, v_1, \ldots, v_{|v|}=v\rangle$ is the
unique geodesic from $\ro $ to $v$. Note that the process is described on the product space $(G,q)^\T$.

A tree-indexed random walk becomes a BRW if the underlying tree is a realization of a Galton--Watson process
with offspring distribution $\mu=(\mu_0,\mu_1,\ldots)$ and mean $m=\sum_k k\mu_k$. For ease of presentation we
will assume that the Galton--Watson process survives almost surely, i.e., $\mu_0=0$, and that $m>1$ in order to
exclude the trivial case $\mu_1=1$. We say the BRW is \emph{recurrent} if $\P(S_v=0 \mbox{ for infinitely many }
v)=1$ and \emph{transient} if $\P(S_v=0 \mbox{ for infinitely many } v)=0$. Here $\P$ does correspond to the
product measure of the Galton--Watson process and the tree-indexed random walk: we pick a realization
$\T(\omega)$ of the Galton--Watson process according to $\mu$ and define the BRW as the tree-indexed random walk
on $\T(\omega)$. Alternatively we could say that the BRW is recurrent if for a.a.~realization $\T(\omega)$ the
tree-indexed random walk is recurrent, i.e., $\P \left( \sum_{n=1}^\infty \sum_{|v|=n} \1\{S_v=\orig\}=\infty
\right)=1$, where $\P$ corresponds to $(G,q)^{\T(\omega)}$. Since we consider BRW on Cayley graphs the
probability $\P(S_v=0 \mbox{ for infinitely many } v)$ is either $0$ or $1$, see~\cite{gantert:04}. Therefore,
the BRW is either recurrent or transient.

We have the following classification due to \cite{gantert:04}:

\begin{thm}\label{thm:rec-tr}
The BRW is transient if and only if $m\le 1/\rho(P)$.
\end{thm}

\begin{rem}\label{rem:defnBRW}
There is the following equivalent description of BRW. At time $0$ we start the process with one particle in
$\orig.$ At time $1$ this particle splits up according to some offspring distribution $\mu$. Then these
offspring particles move (still at time $1$) independently according to $P$. The BRW is now defined inductively:
at each time each particle splits up independently of the others according to $\mu$ and the \emph{new} particles
move then independently according to $P$.
\end{rem}

\section{Dynamical BRW}
Let us introduce the dynamical process. Fix a tree $\T$. For each $v\in \T$, let $\{X_v(t)\}_{t\geq 0}$ be an
independent process that updates its value by an independent sample of $q$ with rate $1$. Formally, consider
i.i.d.~random variables $\{X_v^{(j)}:~v\in\T,~j\in\N\}$ with law $q$, and an independent Poisson process
$\{\psi_v^{(j)}\}_{j\geq 0}$ of rate $1$ for each $v´\in\T$. Define
\begin{equation}\label{eq:xni}
X_v(t):=X_v^{(j)} \mbox{ for } \psi_v^{(j-1)}\leq t < \psi_v^{(j)},\end{equation} where $\psi_v^{(0)}=0$ for
every $n$. The distribution of $(X_v(t))_{v\in\T}$ is $q^\T$ for every $t\geq 0$. Denote $\PP$ the probability
measure on the underlying probability space on which the dynamical BRW process is defined. In the following
$\PP,~\EE$ will always correspond to the dynamical version while $\P,~\E$ describe the non-dynamical process.

Due to Theorem \ref{thm:rec-tr} we have with Fubini's Theorem
\begin{equation*}
\PP \left( \sum_{n=1}^\infty \sum_{|v|=n} \1\{S_v(t)=\orig\}<\infty \mbox{ for Lebesgue-a.e. }t \right)=1
\end{equation*}
if $m\le 1/\rho$ and \begin{equation*} \PP \left( \sum_{n=1}^\infty \sum_{|v|=n} \1\{S_v(t)=\orig\}=\infty
\mbox{ for Lebesgue-a.e. }t \right)=1
\end{equation*} if $m>1/\rho(P)$. The  result of
this note is that there are no exceptional times for transience and recurrence of BRW in the sub-and
supercritical regime. For the critical regime we assume an additional condition on the Cayley graph, but believe
that transience is dynamically stable in general.

\begin{thm}\label{thm:no_except}
We consider a BRW on a Cayley graph $G$ with law $q$ and offspring distribution $\mu$ (whose support excludes
$0$) and  mean $m>1$. Then:
\begin{itemize}
    \item if $m\leq 1/\rho(P)$ and $\sum_{n} n p^{(n)}(\orig,\orig) m^n<\infty$ then \[\PP \left( \sum_{n=1}^\infty \sum_{|v|=n} \1\{S_v(t)=\orig\}<\infty
\mbox{ for all }t \right)=1\]
    \item if $m>1/\rho(P)$ then \[\PP \left( \sum_{n=1}^\infty \sum_{|v|=n} \1\{S_v(t)=\orig\}=\infty
\mbox{ for all }t \right)=1\]
\end{itemize}
\end{thm}

\begin{rem}
Observe that $m<1/\rho(P)$ implies that $\sum_{n} n p^{(n)}(\orig,\orig) m^n<\infty$. The later condition does
in general not hold in the critical case, i.e., $m=1/\rho(P)$: we only have that $\sum_{n} p^{(n)}(\orig,\orig)
m^n<\infty$ for random walks on nonamenable Cayley graphs. The condition depends on the nonexponential type of
the return probabilities; if $p^{(n)}(\orig,\orig)\sim \rho(P)^{n} n^{-\lambda}$ we speak of $n^{-\lambda}$ as
the nonexponential  type of return probabilities. While the simple random walk on the homogeneous tree does not
satisfies the condition, various random walks on free products with $\lambda>2$ are presented in~\cite{CG:09}.
\end{rem}

\begin{rem}\label{rem:moredynam}
In Theorem \ref{thm:no_except} the Galton--Watson tree is fixed for all $t\geq 0$ and only the increments of the
tree-indexed random walk are dependent on $t$. One may use a \emph{dynamical} version of the Galton--Watson
process to define a dynamical BRW. To do this, let $Y(t)$ be a dynamical random variable with values in
$\N=\{1,2,\ldots\}$ defined as in Equation (\ref{eq:xni}) and $(Y_{n,i}(t))_{n,i\geq 1}$ be independent random
variables distributed like $Y(t)$. We define the dynamical Galton--Watson process as
\[ Z_1(t)=1~\mbox{and}~Z_{n+1}(t)=\sum_{i=1}^{Z_n(t)} Y_{n,i}(t)~\mbox{for}~ n\geq 1.\]
Every Galton--Watson process $Z_n(t)$ gives rise to a Galton--Watson tree $\T(t)$ which enables us to define a
dynamical BRW such that the distribution of $(X_v(t))_{v\in\T(t)}$ is $q^{\T(t)}$ for every $t\geq 0$. The proof
of  Theorem \ref{thm:no_except} generalizes to this version of dynamical BRW. Observe hereby that one not only
has to take care about the values of $X_v$'s along certain geodesics but also about the existence of these
geodesics during some time intervals.
\end{rem}

\begin{rem}\label{rem:generalBRW}
A multi-type Galton--Watson process is a generalization of the standard Galton--Watson process where one
distinguishes between different types of particles. The different types are normally encoded by the natural
numbers. The process is described by $Z_n=(Z_n(1), Z_n(2),\ldots)$ where $Z_n(k)$ denotes the number of
particles of type $k$ at time $n$. For $k\geq 1$ let $Y^{(k)}=(Y^{(k)}(1), Y^{(k)}(2),\ldots)$ be random
variables in $\{\N\cup\{0\}\}^\N$. Say we start the process with one particle of type $1$, i.e.,
$Z_1=(1,0,0,\ldots)$, then the multi-type Galton--Watson process is defined inductively by
\[ Z_{n+1}=\sum_{k=1}^\infty \sum_{i=1}^{Z_n(k)} Y^{(k)}_{n,i}~\mbox{for all}~n\geq 1,\]
where $(Y^{(k)}_{n,i})_{n,i\geq 1}$ are independent random variables distributed like $Y^{(k)}$. For ease of
presentation we assume first that the process survives almost sure, i.e., $\sum_{i=1}^\infty Y^{(k)}(i)>1$
almost sure for all $k$. Furthermore, we assume the multi-type Galton--Watson process to be irreducible, i.e.,
for all $k$ there exists some $n\in\N$ such that the probability that $Z_n(k)>0$ is positive.

A BRW on a Cayley graph can be viewed as a multi-type Galton--Watson process: the position $x\in G $ of a
particle is interpreted as its type. In this case the types are encoded by elements of $G$. Recurrence of a BRW
on a Cayley graph means that each site $x$ is visited infinitely many times almost sure. This is equivalent to
say that the corresponding multi-type Galton--Watson process survives locally: $\P(Z_n(x)>0~\mbox{for infinitely
many}~n)=1$. If we leave the homogeneous setting of BRW on Cayley graphs and consider BRW on a general graph the
probability that infinitely many particles return to the starting position may be strictly between $0$ and $1$,
compare with \cite{gantert:04}. We say the BRW is transient if the latter probability is $0$ and recurrent if it
is positive. The corresponding phenomenon occurs for multi-type Galton--Watson processes too. We say there is
local survival if $\P(Z_n(k)>0~\mbox{for infinitely many}~n)>0$. Due to the irreducibility of the process the
latter probability is either positive for all $k$ or equal to $0$ for all $k$ and we can speak of local
extinction if $\P(Z_n(k)>0~\mbox{for infinitely many}~n)=0$ for all (some) $k$, compare with \cite{gantert:09}.

It is straightforward to define a dynamical multi-type Galton--Watson process analogously to the dynamical
Galton--Watson process defined in Remark \ref{rem:moredynam}. One might wonder if there are exceptional times
for local survival and local extinction. Now, let us drop the assumption that the process survives almost
surely. We speak of global survival if $Z_n>0$ for infinitely many $n$ with positive probability and of global
extinction otherwise. The treatment of global survival is in general more difficult and even more subtle since
\emph{critical} processes may survive or die out, compare with \cite{BZ:08}. Therefore the study of exceptional
times for global survival/extinction is one of the next steps to go.
\end{rem}

\begin{rem}
Let us consider a transient random walk $S_n=\sum_{i=1}^n X_i$ on $\Z$ (or $\R$) with $\E[X_i]>0.$ We assume
that  there exists a rate function $I(\cdot)$ satisfying
\[ -I(a)=\lim_{n\to\infty} \frac1n \log \P(S_n\leq an) \mbox{ for
} a\leq \E[X_i].\] Denote by $m_n$ the minimal position of a particle at time $n$; $m_n=\min_{|v|=n} S_v$. There
is the classical result that $\lim_{n\to\infty} \frac{m_n}{n}=\inf\{s:~ I(s)\leq \log m\}$. Combining the proof
of Theorem 18.3 in \cite{peres99} with the ideas of the proof of Theorem \ref{thm:no_except} one can see that
there are no exceptional times for the (linear) speed, i.e., $\lim_{n\to\infty} \frac{m_n(t)}{n}=\inf\{s:~
I(s)\leq \log m\}$ for all $t$. Furthermore, in the critical case $m=1/\rho(P)$ we have that $m_n/n\to 0$ but
$m_n\to \infty$ as $n\to \infty$. The \emph{second order} behaviour is more subtle: while for a wide range of
BRW $m_n/\log n$ converges in probability it does in general not converge almost surely. We refer to \cite{hu08}
for more details and references on recent results. In this respect the study of exceptional times for the second
order behaviour is of interest.
\end{rem}

\section{Proof of Theorem \ref{thm:no_except}}

\noindent\emph{Transience is dynamically stable if $\sum n p^{(n)}(\orig, \orig) m^n<\infty$:} It is convenient
to define an auxiliary random variable $\tau$, which is exponentially distributed with mean $1$ and independent
of $(X_{v}(t))_{v\in \T, t\geq 0}$, see Section~3 in \cite{MR1959784}. Define
\[ \ZZ_n:= \int_0^\tau \sum_{|v|=n} \1\{S_v(t)=\orig\} dt.\]
By Fubini's Theorem   we have $\EE[\ZZ_n]=m^n p^{(n)}(\orig,\orig)$. We follow the line of proof of Lemma 5.6 in
\cite{MR1959784}. In what follows we only consider those $n\in\N$ such that $\P(S_n=e)>0$. We have for $n\geq 1$

\[\PP(\ZZ_n>0)=\frac{\EE[\ZZ_n]}{\EE[\ZZ_n|\ZZ_n>0]}. \]
Let $\sigma:=\inf\{ t\geq 0: S_v(t)=\orig\mbox{ for some } v\in \T_n\}.$  There might be several $v$ such that
$S_v(\sigma)=\orig$. In order to choose one of them we use an enumeration of the vertices of $\T_n$. This
enumeration is chosen once for $t=0$ and remains fixed for  $t\geq 0$. Define the random set $R=\{v\in\T_n:
S_v(\sigma)=\orig\}$ and let $\tilde v$ be the smallest (in the above enumeration) element of $R$. Furthermore,
denote $\langle \ro, \tilde v_1,\ldots,\tilde v_n=\tilde  v\rangle$ the geodesic from $\ro$ to $\tilde v$. For
any geodesic $\langle\ro, v \rangle=\langle \ro,  v_1,\ldots, v_n=v\rangle$ we define the event
\[ A(\langle\ro, v \rangle)=\{ X_{v_k},~1\le k\le n,~\mbox{do not change their values during}~[\sigma,\sigma+1/n]\}.\]
Since the event that $X_{v_k}$ changes its value during $[\sigma,\sigma+1/n]$ is independent of its value at
time $\sigma$ we have that $\PP(A(\langle\ro, v \rangle)\mid \sigma<\infty)=1/e$.

Conditioned on the event $\{\ZZ_n>0\}$ we have $\sigma\in[0,\tau)$. By the strong Markov property and the
\emph{memoryless} property of $\tau$, we have
\begin{eqnarray}\label{eq:dyn}
& & \PP( \tau>\sigma +1/n, \exists v:~ S_v(t)=\orig~\forall t\in[\sigma, \sigma+1/n]|\ZZ_n>0)\cr  & &~\geq
\PP(\tau>\sigma +1/n, A(\langle \ro, \tilde v)\rangle)|\ZZ_n>0)= \left(\frac1e\right)^{1/n} \frac1e.
\end{eqnarray}
If the above event occurs, then $\ZZ_n\geq 1/n$. Hence
\[\PP(\ZZ_n\geq 1/n|\ZZ_n>0)\ge 1/e^2  \mbox{ and }  \EE[\ZZ_n|\ZZ_n>0] \geq \frac1{ne^2}.\]
Eventually,
\[\PP(\ZZ_n>0)\leq e^2 n  p^{(n)}(\orig,\orig)m^n\quad \forall n,\]
and hence $\sum_n \PP(\ZZ_n>0)<\infty$. By the lemma of Borel--Cantelli there are no times $t$ such that
$\sum_{|v|=n} \1\{S_v(t)=\orig\}>0$ for infinitely many $n$ and therefore
\[ \PP \big(~\exists t:~ \sum_n \sum_{|v|=n}
\1\{S_v(t)=\orig\}=\infty\big)=0.\]

\noindent\emph{Recurrence is dynamically stable:} Let us first consider the non-dynamical BRW. Since
$m>1/\rho(P)$ we have with equation (\ref{eq:rho:prop}) that there exists some $k\in\N$ such that
$p^{(k)}(\orig,\orig) m^k >1$. It's rather standard, e.g. compare with proof of Theorem 18.3 in \cite{peres99},
to define the following embedded process $(\xi_n)_{n\ge 1}$ where we observe the process only at times
$ik,~i\in\N$, and kill all particles that are not in $\orig$ at these times. Then $\xi_n$ describes the number
of particles at $\orig$ at time $nk$. Let us give a more formal definition of this process. Let $\T^v$ be the
induced subgraph of $\T$ consisting of all descendants of $v$, or in other words, of all vertices whose geodesic
to the root goes through $v$. Furthermore, let $\T^v_k=\{w\in \T^v:~ d(v,w)=k\}$ be the vertices at level $k$ of
$\T^v$ and denote by $v^{-k}$ the $k$-th predecessor of $v$, i.e., $v\in \T^{v^{-k}}$ and $d(v^{-k},v)=k$.
Define
\[Y_v=\sum_{w\in T^v_k} \1\{S_w=S_v\}.\] Let $H_1=\{\ro\}$ and define inductively for $n\geq 1$
\[H_{n+1}=\{v:~|v|=nk,~ S_v=\orig,~S_{v^{-k}}\in H_{n}\}.\]
Then $\xi_n=|H_n|$ but can also be written as
\[ \xi_{n+1}=\sum_{v\in H_{n}} Y_v.\]
Since the $((Y_v)_{v\in\T_{ik}})_{i\geq 1}$ are i.i.d. random variables the process $\xi_n$ has the same law as
the Galton--Watson process $(Z_n)_{n\ge 1}$ defined through $Z_1=1$ and $Z_{n+1}=\sum_{i=1}^{Z_{n}} Y_{n,i},$
where $(Y_{n,i}(_{n,i\geq 1}$ are i.i.d. random variables distributed like $Y_v$. Since
$\E[Y_v]=p^{(k)}(\orig,\orig)m^{k}>1$ the process $Z_n$ is a supercritical Galton--Watson process and hence
survives with positive probability and so does $\xi_n$.

The dynamical version of the BRW induces a dynamical version of $\xi_n$: Let $H_1(t)=\{\ro\}$ and define
inductively
\[H_{n+1}(t)=\{v:~|v|=nk,~ S_v(t)=\orig,~S_{v^{-k}}(t)\in
H_{n}(t)\}.\] Let $\xi_n(t)=| H_n(t)|$ which can also be written as
\[ \xi_{n+1}(t)=\sum_{v\in H_{n}(t)} Y_v(t)\] where \[Y_v(t)=\sum_{w\in T_v^k} \1\{S_w(t)=S_v(t)\}.\]
Again, with $(Y_{n,i}(t))_{n,i\geq 1}$ i.i.d. like $Y_v(t)$ the Galton--Watson process defined by $Z_1(t)=1$ and
\[Z_{n+1}(t)=\sum_{i=1}^{Z_{n}(t)} Y_{n,i}(t),\]
has the same law as $\xi_n(t)$.

As before $\P$ denotes the probability measure for the non-dynamical process while $\PP$ describes the dynamical
version. Analogously to the study of the noncritical cases in \cite{haggstrom97} let
\[\inf_{[a,b]} H_n = \bigcap_{t\in[a,b]} H_n(t)~\mbox{and}~\inf_{[a,b]}\xi_n=|\inf_{[a,b]} H_n|.\]
With
\[\inf_{[a,b]}Y_{v}=\sum_{w\in\T^v_k} \1\{S_w(t)=S_v(t)~\forall t\in[a,b]\}\]
we can write
\[ \inf_{[a,b]}\xi_n = \sum \inf_{[a,b]}Y_{v},~\mbox{ where the sum is over}~v\in
\inf_{[a,b]}H_{n-1}.\]

As above, the process $\inf_{[a,b]} \xi_n$ is a Galton--Watson process with mean $\EE[ \inf_{[a,b]}Y_{v}]$. In
order to show that there exists some $\eps>0$ such that $\EE[ \inf_{[0,\eps]}Y_{v}]>1$ we proceed analogously to
the arguments around Equation (\ref{eq:dyn}). Let $R_v=\{w\in T_k^v:~ S_w(0)=S_v(0)\}.$ Condition on the fact
that $|R_v|\geq l$ let $w_1,\ldots,w_l$ be the $l$ smallest (according to some enumeration of the vertices)
elements of $R_v$. Observe
\begin{eqnarray*}
\PP(\inf_{[0,\eps]} Y_v \ge l)&\geq&  \PP(|R_v|\geq l) \PP(S_{w_i}(t)=S_{w_i}(0)~ \forall t\in[0,\eps], 1\le i
\le l\mid |R_v|\geq l)\cr &\geq & \P( Y_v\ge l) \left(\frac1e\right)^{\eps k l}
\end{eqnarray*}
where the last inequality comes from  the consideration of the \emph{worst case} where the geodesics $\langle
\ro, w_i\rangle$ are disjoint (except of $\ro$).

Eventually, since $\E[Y_v]>1$ we can choose $\eps>0$ such that $\EE [\inf_{[0,\eps]} Y_{n,i}] >1.$ Hence, the
Galton--Watson process $\inf_{[0,\eps]} \xi_n$ is supercritical and there are no exceptional times in the
interval $[0,\eps]$. Repeating the arguments  for the intervals $[k\eps,(k+1)\eps]$ and using countable
additivity concludes the proof for $m>1$. \vspace{1cm}

\noindent \textbf{Acknowledgment:} Thanks to Itai Benjamini for suggesting this problem and introducing me to
the concept of dynamical sensitivity. The author wishes to thank a referee for helpful suggestions and for
bringing to his attention several inaccuracies.

\begin{small}
\addcontentsline{toc}{chapter}{Bibliography}
\bibliography{bib}

\begin{thebibliography}{10}

\bibitem{MR1959784}
Itai Benjamini, Olle H{\"a}ggstr{\"o}m, Yuval Peres, and Jeffrey~E. Steif.
\newblock Which properties of a random sequence are dynamically sensitive?
\newblock {\em Ann. Probab.}, 31(1):1--34, 2003.

\bibitem{benjamini:94}
Itai Benjamini and Yuval Peres.
\newblock Markov chains indexed by trees.
\newblock {\em Ann. Probab.}, 22(1):219--243, 1994.

\bibitem{BZ:08}
Daniela Bertacchi and Fabio Zucca.
\newblock Characterization of the critical values of branching random walks on
  weighted graphs through infinite-type branching processes, arXiv:0804.0224v1.

\bibitem{CG:09}
Elisabetta Candellero and Lorenz~A. Gilch.
\newblock Phase transitions for random walk asymptotics on free products of
  groups, arXiv:0909.1893v2.

\bibitem{gantert:04}
Nina Gantert and Sebastian M\"uller.
\newblock The critical branching {M}arkov chain is transient.
\newblock {\em Markov Process. and Rel. Fields.}, 12:805--814, 2006.

\bibitem{gantert:09}
Nina Gantert, Sebastian M\"uller, Serguei Popov, and Marina Vachkovskaia.
\newblock Survival of branching random walks in random environment.
\newblock {\em to appear in Journal of Theoretical Probability}, DOI
  10.1007/s10959-009-0227-5, 2009+.

\bibitem{haggstrom97}
Olle H{\"a}ggstr{\"o}m, Yuval Peres, and Jeffrey~E. Steif.
\newblock Dynamical percolation.
\newblock {\em Ann. Inst. H. Poincar\'e Probab. Statist.}, 33(4):497--528,
  1997.

\bibitem{MR2235173}
Christopher Hoffman.
\newblock Recurrence of simple random walk on {${\Bbb Z}\sp 2$} is dynamically
  sensitive.
\newblock {\em ALEA Lat. Am. J. Probab. Math. Stat.}, 1:35--45 (electronic),
  2006.

\bibitem{hu08}
Yueyun Hu and Zhan Shi.
\newblock Minimal position and critical martingale convergence in branching
  random walks, and directed polymers on disordered trees.
\newblock {\em {A}nnals of {P}robability}, to appear in.

\bibitem{Khoshnevisan05}
Davar Khoshnevisan, David~A. Levin, and Pedro~J. M{\'e}ndez-Hern{\'a}ndez.
\newblock On dynamical {G}aussian random walks.
\newblock {\em Ann. Probab.}, 33(4):1452--1478, 2005.

\bibitem{Khoshnevisan06}
Davar Khoshnevisan, David~A. Levin, and Pedro~J. M{\'e}ndez-Hern{\'a}ndez.
\newblock Exceptional times and invariance for dynamical random walks.
\newblock {\em Probab. Theory Related Fields}, 134(3):383--416, 2006.

\bibitem{peres99}
Yuval Peres.
\newblock Probability on trees: an introductory climb.
\newblock In {\em Lectures on probability theory and statistics (Saint-Flour,
  1997)}, volume 1717 of {\em Lecture Notes in Math.}, pages 193--280.
  Springer, Berlin, 1999.

\bibitem{steif09}
Jeffrey~E. Steif.
\newblock A survey on dynamical percolation.
\newblock {\em arXiv:0901.4760}, 2008++.

\bibitem{woess}
Wolfgang Woess.
\newblock {\em Random walks on infinite graphs and groups}, volume 138 of {\em
  Cambridge Tracts in Mathematics}.
\newblock Cambridge University Press, Cambridge, 2000.

\end{thebibliography}
\end{small}
\end{document}